\documentclass{article}

\usepackage{tabularx}
\usepackage{graphicx}
\usepackage{times}

\usepackage{authblk}

\usepackage{amsmath}
\newcommand\numberthis{\addtocounter{equation}{1}\tag{\theequation}}

\usepackage{hyperref}

\usepackage{amssymb}
\usepackage{xcolor}

\title{Understanding the impact of numerical solvers on inference for differential equation models}

\author[1]{Richard Creswell}
\author[1]{Katherine M.\ Shepherd}
\author[2]{Ben Lambert}
\author[3]{Gary R.\ Mirams}
\author[4]{Chon Lok Lei}
\author[5]{Simon Tavener}
\author[1]{Martin Robinson}
\author[1]{David J.\ Gavaghan}
\affil[1]{Department of Computer Science, University of Oxford}
\affil[2]{Department of Mathematics, University of Exeter}
\affil[3]{School of Mathematical Sciences, University of Nottingham}
\affil[4]{Institute of Translational Medicine, Faculty of Health Sciences, University of Macau; Department of Biomedical Sciences, Faculty of Health Sciences, University of Macau}
\affil[5]{Department of Mathematics, Colorado State University}
\date{}

\begin{document}
\maketitle

\begin{abstract}\noindent
Most ordinary differential equation (ODE) models used to describe biological or physical systems must be solved approximately using numerical methods. Perniciously, even those solvers which seem sufficiently accurate for the \emph{forward problem}, i.e., for obtaining an accurate simulation, may not be sufficiently accurate for the \emph{inverse problem}, i.e., for inferring the model parameters from data. We show that for both fixed step and adaptive step ODE solvers, solving the forward problem with insufficient accuracy can distort likelihood surfaces, which may become jagged, causing inference algorithms to get stuck in local ``phantom'' optima. We demonstrate that biases in inference arising from numerical approximation of ODEs are potentially most severe in systems involving low noise and rapid nonlinear dynamics. We reanalyze an ODE changepoint model previously fit to the COVID-19 outbreak in Germany and show the effect of the step size on simulation and inference results. We then fit a more complicated rainfall-runoff model to hydrological data and illustrate the importance of tuning solver tolerances to avoid distorted likelihood surfaces. Our results indicate that when performing inference for ODE model parameters, adaptive step size solver tolerances must be set cautiously and likelihood surfaces should be inspected for characteristic signs of numerical issues.
\end{abstract}

\textbf{Keywords:} ordinary differential equations, inference, Bayesian statistics, truncation error,  compartmental models, hydrological modelling

\section{Introduction}

Many scientific phenomena involve time-varying signals or outputs. These phenomena are often believed to obey some parametric model, whose parameter values are \emph{a priori} unknown but can be inferred from observed data. In this paper, we present and analyze the key challenges that arise in parameter inference when models involve ordinary differential equations (ODEs). ODEs are used throughout the biological and physical sciences to express dynamic processes; a few examples amongst the myriad of their application areas include epidemiology~\cite{van2022learning}, hydrology~\cite{kavetski2003semidistributed}, cardiac electrophysiology~\cite{whittaker2020calibration}, and population dynamics~\cite{shertzer2002predator}.

In general, the differential equations used in scientific applications cannot be solved analytically. However, a wide range of computational methods have been developed to obtain a numerical approximation of their solutions. (Solving the differential equation for a particular value of the parameters is known as the \emph{forward problem}.) While numerical algorithms to solve the forward problem introduce error, the properties of this error are generally well understood and can be controlled. In solvers using a fixed time step (discussed in \S{}\ref{fixed_step}), the error can be reduced by decreasing the size of the time step~\cite{gautschi1997numerical}. In solvers in which the time step is set adaptively (discussed in \S{}\ref{adaptive}), the error is typically controlled through user-specified relative and absolute tolerances on the local truncation error (the error in the solution introduced by a single time step of the solver)~\cite{dormand1980family}.

Our focus in this paper is on the interplay between the numerical approximations inherent in the forward problem and the \emph{inverse problem}, which consists of learning values of the parameters that are compatible with an observed time series. Widely used approaches to the inverse problem include optimization of an objective function which measures the quality of fit between the model and the data (e.g., maximum likelihood), or Bayesian approaches which generate samples from the posterior distribution of the parameters (e.g., Markov chain Monte Carlo (MCMC)). These approaches to the inverse problem require the forward problem to be solved at multiple different parameter values. The errors in each numerical solution of the forward problem, even when individually small, are liable to introduce significant bias to inference results.


The rest of this paper is organised as follows. In \S{}\ref{effects}, we present the widely used independent and identically distributed Gaussian noise log-likelihood function for fitting ODE models and derive a bound on the error in this log-likelihood arising from the use of an approximate solution to the ODE. On the basis of this bound, and results presented subsequently, we argue that the biases in inference results arising from numerical solvers are likely to be most severe in systems which have low noise and rapid nonlinear dynamics. In \S{}\ref{sec:3}, we study two broad classes of ODE solvers: those involving a fixed time step, and those involving a time step set adaptively to control the error on the solution. For both classes of solver, we illustrate the effects that solver inaccuracy may have on inference, and illustrate this using synthetic data. Additionally, in \S{}S3, we study how smoothing approximations can reduce the influence of numerical error on computation of the likelihood. Finally, in \S{}\ref{sir} and \S{}\ref{hydrology} we consider inference of ODE models for real data series. in \S{}\ref{sir} we reanalyze an ODE model of disease transmission fit to the COVID-19 outbreak in Germany and show that, when using a solver with a fixed time step, the choice of time step can alter inference and simulation results, and in \S{}\ref{hydrology} we fit a rainfall-runoff model to hydrological data to illustrate the pitfalls of performing parameter inference using an adaptive step size solver with insufficient local tolerances.

\section{Effects of numerical error on computation of the log-likelihood}  \label{effects}

\subsection{Log-likelihood function for an ODE model}

We assume that time series data $\{y_i\}_{i=1}^N$; $y_i \in \mathbb{R}^n$ are measured at time points $\{t_i\}_{i=1}^N$. These data are believed to obey some function $g: \mathbb{R}^l \rightarrow \mathbb{R}^n$ of $x(t; \theta) \in \mathbb{R}^l$, where $x$ is the solution to an ordinary differential equation:
\begin{align}  \label{eq:ode}
    \frac{dx}{dt} &= h(t, x, \theta); \\
    x(0) &= x_0,
\end{align}
for some function $h$ which is informed by the relevant scientific theory and parameterized by the (potentially unknown) values $\theta \in \mathbb{R}^m$. (In some inference problems, the initial value $x_0$ is also unknown and inferred from the data.) Eq.\ \eqref{eq:ode} has been written as a \emph{first order} equation (i.e., only involving the first order derivative of $x$); higher order differential equations may be rewritten as systems of first order equations.

Deterministic forward models never fully explain the variation in real observations. To make the forward model a feasible explanation of the data, an additional stochastic component representing unmodelled elements (for example, processes involved in the measurement of the signal) is included. This stochastic component is often included in additive form, yielding the proposed data generating process:
\begin{equation}  \label{eq:noise}
y_i = g(x(t_i; \theta)) + \epsilon_i,
\end{equation}
where $\epsilon_i$ is a mean-zero random variable specifying the noise process. Many choices for $\epsilon_i$ are possible, and the assumed distribution of this variable should be chosen carefully, as misspecification of $\epsilon_i$ may cause inaccurate inference results \\\cite{lambert2022autocorrelated}. A standard choice, which we use for synthetic data generation in this paper, is the independent and identically distributed (IID) Gaussian:
\begin{equation}  \label{eq:iid_gaussian}
\epsilon_i \overset{\text{IID}}{\sim} N(0, \sigma),
\end{equation}
with the parameter $\sigma$ representing the standard deviation of the noise process, which is also inferred from the data together with the model parameters $\theta$. Choosing a particular noise process enables the joint probability of the observations to be expressed as a function of the parameter values $\theta$---this is known as a likelihood. For an IID Gaussian noise process (eq.~\eqref{eq:iid_gaussian}), the log-likelihood for time series data $\{y_i\}_{i=1}^N$ takes the form:
\begin{equation}  \label{eq:iid_gaussian_ll}
\log p(y_1, \dots, y_N|\theta,\sigma) = -\frac{N}{2} \log(2\pi) - \frac{N}{2} \log(\sigma^2) - \frac{1}{2\sigma^2} \sum_{i=1}^N (y_i - g(x(t_i; \theta)))^2.
\end{equation}

The likelihood expresses the quality of the fit between the model output and the data, with higher values of the likelihood indicating a superior fit. Thus, the values of $\theta$ most compatible with data can be found by maximising the likelihood with respect to $\theta$ (i.e., the method of maximum likelihood). Alternatively, the likelihood can be used together with a prior distribution on the parameters ($p(\theta)$) to infer the posterior distribution according to Bayes' theorem:
\begin{equation*}
p(\theta | y_1, \dots, y_N) = \frac{p(y_1, \dots, y_N|\theta) p(\theta)}{p(y)}.
\end{equation*}
In this paper, we consider the typical case where the posterior cannot easily be expressed in closed form but can be approximated using Markov chain Monte Carlo (MCMC) sampling methods~\cite{gelman2013bayesian}.

\subsection{Error in the log-likelihood arising from approximation of the forward solution}

The data are assumed to obey the IID Gaussian log-likelihood, eq.\ \eqref{eq:iid_gaussian_ll}.
We assume that $x(t_i; \theta)$ is the \emph{true} solution to the ODE at time point $t_i$, which is unavailable and approximated by $\hat{x}_i$. The deviation between $x(t_i; \theta)$ and $\hat{x}_i$ at any time point is given by the global truncation error, $e(t_i)$: $$e(t_i) = x(t_i; \theta) - \hat{x}_i.$$
In general, $e(t_i)$ is unknown, although, for particular numerical solvers, its magnitude may be bounded by some function of the step size or some other quantity which can be used to tune the accuracy of the solver.

The log-likelihood available to the inference algorithm takes the same form as eq.\ \eqref{eq:iid_gaussian_ll}, but computed using the numerical approximation $\hat{x}_i$ instead of $x(t_i; \theta)$. For brevity, we denote the accurate log-likelihood by $\mathcal{L}$, and we denote the log-likelihood computed using a numerical approximation by $\mathcal{L}'$, which is given by
\begin{equation} 
\mathcal{L}' = -\frac{N}{2} \log(2\pi) - \frac{N}{2} \log(\sigma^2) - \frac{1}{2\sigma^2} \sum_{i=1}^N (y_i - g(\hat{x}_i))^2.
\end{equation}
Assuming Lipschitz continuity of the observation function $g$ with Lipschitz constant $K$, we prove in \S{}S1 that the difference between $\mathcal{L}'$ and $\mathcal{L}$ is bounded according to:
\begin{equation} \label{eq:ll_error}
    |\mathcal{L}-\mathcal{L}'| \leq  \sum_{i=1}^N \left( \frac{K^2}{2 \sigma^2} |e(t_i)|^2 + \frac{K}{\sigma^2} |e(t_i)| |y_i - g(x(t_i; \theta))|  \right).
\end{equation}
We observe an inverse relationship between $\sigma$ and the bound of $|\mathcal{L}-\mathcal{L}'|$ when $e(t_i)$ is held constant. Thus, when a solver is tuned to yield a particular global truncation error $e(t_i)$, we expect the absolute bias in the log-likelihood to be more severe at smaller values of $\sigma$. Eq.\ \eqref{eq:ll_error} also indicates that more severe biases may be expected when $N$ is larger---i.e., there are more observations. Additionally, at a fixed level of $\sigma$, we expect the bias in the log-likelihood to decrease as the global truncation error is decreased. In \S{}S2, we additionally derive the distribution of the error in the likelihood, and show that $\mathbb{E}[\mathcal{L}-\mathcal{L}']>0$: the numerical approximation of the likelihood will, on average, underestimate the true likelihood.


\section{Effects of ODE solvers on inference}  \label{sec:3}

To study the interplay between ODE solvers and inference, we introduce the following differential equation problem which describes an oscillatory system with damping and forcing:
\begin{equation*}  \label{eq:oscillator}
m \frac{d^2x}{dt^2} + c \frac{dx}{dt} + k x = F(t).
\end{equation*}
The model has three parameters which will be treated as unknown: $(m, c, k)$. In classical mechanics, these represent the mass, damping coefficient, and spring constant respectively. $F(t)$ represents the forcing function or stimulus, and in this paper takes a variety of forms throughout our results. This damped and forced oscillator is described by a second order differential equation; to apply ODE solvers straightforwardly, we rewrite it as a first order differential equation of two state variables:
\begin{equation}  \label{eq:oscillator_2}
\frac{d}{dt} \begin{pmatrix} x \\ \dot{x} \end{pmatrix} = \begin{pmatrix}
\dot{x} \\ \frac{F(t)}{m} - \frac{c}{m} \dot{x} -\frac{k}{m} x
\end{pmatrix},
\end{equation}
where $\dot{x}=dx/dt$.

\subsection{Fixed step and adaptive step ODE solvers}

A wide range of numerical algorithms have been developed to obtain approximate solutions to initial value problems (IVPs) of the form given in eq.\ \eqref{eq:ode}. These algorithms typically work by computing an approximate solution on a grid of time points (in general, distinct from the time points where the data are located) and then using an interpolation algorithm to obtain the solution at intermediate time points.

Most simply, the grid of solver time points can be prespecified in advance (we refer to such methods as fixed time step solvers). However, in general, it is inefficient to use the same time step throughout the entire time range on which the ODE is being solved, particularly when solved repeatedly over a range of parameters. Solvers can employ large time steps in regions where the solution and its gradients change gradually without causing much error in the solution; however, in regions where the derivative changes rapidly, small time steps are required to maintain a low error. This motivated the development of ODE solvers which adjust the step size throughout the time domain over which the ODE is solved. While fixed step solvers are still commonly used, adaptive step solvers are standard in high performance computing and are widely implemented in software libraries for ODE solving.

When using an adaptive step size solver, the user does not specify a step size, but rather a local error tolerance. The algorithm then selects a time-varying sequence of step sizes such that the local error in the solution falls below the specified tolerance. The total number of time steps used by the solver thus depends on the selected tolerance and the properties of the solution. Typically, an interpolation scheme is then used to obtain the solution at intermediate time values. Tolerances can be expressed either as an absolute value or relative to the magnitude of the solution. In many implementations, both are available to the user: for example, the SciPy library allows the user to specify both an absolute tolerance \texttt{atol} and a relative tolerance \texttt{rtol}, and chooses step sizes such that the magnitude of the local truncation error on the solution $x$ does not exceed $\texttt{atol} +\texttt{rtol} | x |$~\cite{2020SciPy-NMeth}. For the results presented in this paper, we fix $\texttt{atol}$ to a value of $10^{-9}$ and tune $\texttt{rtol}$ to control the accuracy of the solver. Adaptive step sizes have been implemented for a wide variety of ODE solver algorithms. Here, we focus on Runge-Kutta methods of the form RK$p(q)$, which use the $q$th order method to estimate the error (and thus control the time step), while making the actual steps using the $p$th order method~\cite{dormand1980family}. Runge-Kutta methods are not described in detail here for brevity---they are widely used and details can be found in many standard texts (for example,~\cite{gautschi1997numerical}). We rely on the SciPy adaptive time step Runge-Kutta implementation, which employs a quartic interpolation polynomial for RK5(4) and a cubic Hermite interpolation polynomial for RK3(2) \cite{2020SciPy-NMeth}.

\subsubsection{Typical log-likelihood surface shapes}
We now consider the influence of the two numerical solution methods for parameter inference. Because fixed step solvers use the same grid throughout parameter space, while adaptive step solvers may employ different grids at different parameter values, these two classes of solvers differ in the characteristics of the error that they may introduce into the likelihood function. 

\begin{figure}
\centering
\includegraphics[width=.75\textwidth]{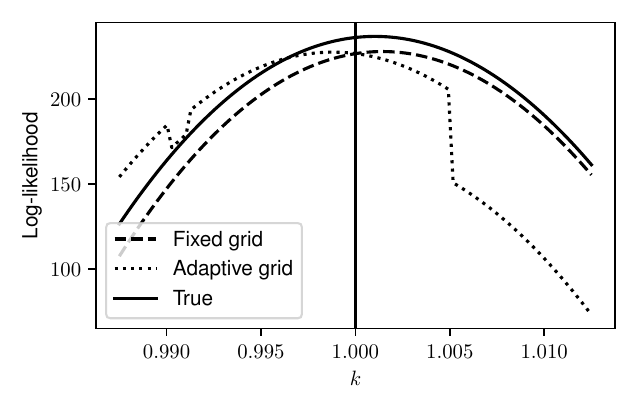}
\caption{\textbf{Comparison of log-likelihood surfaces calculated using fixed step and adaptive step solvers.} Log-likelihood for the parameter $k$ calculated from data generated from the oscillator model eq.\ \eqref{eq:oscillator_2}, with all other parameters held at their true values. The log-likelihood was calculated from eq.~\eqref{eq:iid_gaussian_ll} using an adaptive step RK5(4) solver with relative tolerance set to $10^{-8}$ (True), a Forward Euler solver with a fixed time step $\Delta t=0.01$, and an adaptive step RK5(4) solver with tolerance tuned such that at the true parameter values (vertical line) it introduces the same magnitude of error in the log-likelihood as the fixed step Forward Euler solver (corresponding to a relative tolerance of $0.00944$).}
\label{fig:compare}
\end{figure}

We illustrate this by computing the likelihood surface for the $k$ parameter in the oscillator problem, eq.\ \eqref{eq:oscillator_2}. $75$ evenly spaced data points were generated from and including $t=0$ to $t=50$ from the model with an accurate solver (the RK5(4) solver with relative tolerance set to $10^{-8}$), using true parameter values $k=1$, $c=0.2$, $m=1$, initial conditions of $x(t=0)=0$, $\dot{x}(t=0)=0$, and $$F(t)=\begin{cases}1, & t<25, \\ 0.9, & t\geq 25.\end{cases}$$ Then IID Gaussian noise was added to the solution at each of the sampled locations with $\sigma=0.01$. Holding all other parameters fixed at their true values, the log-likelihood was calculated for a range of values of $k$, using three different ODE solvers. First, the RK5(4) solver with relative tolerance set to $10^{-8}$ was used to compute the accurate (`True') likelihood. Next, the Forward Euler solver with a fixed time step of $\Delta t = 0.01$ was used. Finally, we used the RK5(4) solver, but with its relative tolerance tuned so the observed magnitude in the error in the log-likelihood at the true parameter values was equal to that produced by the Forward Euler solver (for this problem, this resulted in relative tolerance tuned to $0.00944$). These results are shown in Figure \ref{fig:compare}.

At the true parameter value, both solvers result in a slight underestimation of the log-likelihood. Across the parameter range considered, the fixed time step solver results in a log-likelihood which is shifted relative to the true one, but retains the smooth, unimodal shape. However, the adaptive step solver results in a log-likelihood surface which in addition to being shifted exhibits jagged, discontinuous fluctuations. In the remainder of \S{}\ref{sec:3}, we examine these two phenomena in more detail. 

\subsection{Fixed time step solvers} \label{fixed_step}

\subsubsection{Forward Euler solver} 


One of the simplest numerical solvers for ordinary differential equations is the Forward Euler method with a uniform step size $\Delta t$. This solver is easily implemented and thus has achieved wide usage despite its simplicity and typically mediocre performance.

Forward Euler has been used for inference in some recent high-profile epidemiological research where $\Delta t$ was set to a value comparable to the time scale of the behavior of the system (e.g.,~\cite{birrell2021real, dehning2020inferring}). Whether these applications are representative of the use of Forward Euler more generally is unclear, but our results in \S{}\ref{sir} indicate that such choices of $\Delta t$ may alter both forward model solutions and parameter inference results.

\subsubsection{Inference for the damped, driven oscillator using Forward Euler}  \label{fixed_ode}

We now exemplify the impact of using Forward Euler with insufficiently small time steps on inference by using synthetic noisy data generated from the (accurate) solution of eq.\ \eqref{eq:oscillator_2}. $25$ evenly spaced data points were generated from and including $t=0$ to $t=5$ from the model with an accurate solver (the RK5(4) solver with relative tolerance set to $10^{-8}$), using true parameter values $k=1$, $c=0.2$, $m=1$, an initial condition of $x(t=0)=0$, $\dot{x}(t=0)=0$, and $F(t)=1$. Then, IID Gaussian noise was added to the solution at each of the sampled locations with $\sigma=0.1$. Holding all other parameters fixed at their true values, the log-likelihood was calculated for a range of values of $k$, using the Forward Euler solver with various time steps.

\begin{figure}
\centering
\includegraphics[width=\textwidth]{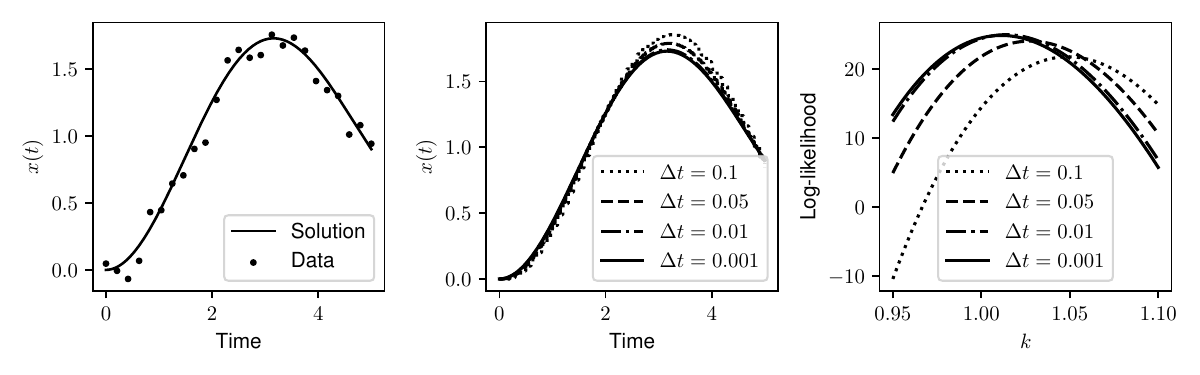}
\caption{\textbf{Damped oscillator inference using Forward Euler.} (Left) Synthetic data for the damped driven oscillator. The curved line indicates the accurate solution to the ODE with these parameters, while the points indicate the noisy data. (Center) Solution for oscillator computed using a Forward Euler solver with four different choices for the time step $\Delta t$. (Right) Log-likelihood for the parameter $k$ calculated from the noisy data, with all other parameters held at their true values. The log-likelihood was calculated from eq.~\eqref{eq:iid_gaussian_ll} using a Forward Euler solver with four different choices for the time step $\Delta t$.}
\label{fig:fixed_step}
\end{figure}

Figure~\ref{fig:fixed_step} shows the impact of using Forward Euler on the likelihood surface. The results show the typical effect of a fixed step solver with insufficiently small time steps: the likelihood surface maintains a smooth shape, but it is shifted relative to its true location. The longest time step considered in this study, $\Delta t = 0.1$, causes substantial inaccuracy in the likelihood even though $\Delta t=0.1$ is small compared to the time scale of the dynamics of the system and the system with $F(t)=1$ contains no discontinuities or other challenging features.

As the step size is refined, the log-likelihood curves converge. This suggests a diagnostic technique which could be incorporated into inference algorithms: once the optimal parameter values have been determined, the log-likelihood should be evaluated at those parameter values with the step size on the solver slightly adjusted; if the solver is sufficiently accurate, the value of the log-likelihood should not be a strong function of the step size.

\subsection{Adaptive step size solvers}  \label{adaptive}

Adaptive step size solvers enable increased efficiency in obtaining accurate solutions to ODEs. However, when used in inference problems, they can convert a smooth likelihood surface into a rough one, characterized by rapid (and entirely phantom) changes in the likelihood which interfere with inference algorithms. These inaccuracies in the likelihood can be observed even at tolerances in the solution error where further refinements do not visibly influence the solution. For example, in cardiac electrophysiology, jagged parameter likelihoods have been observed with adaptive step size ODE solvers with tolerances as low $10^{-7}$~\cite{johnstone2018uncertainty, gary_blog}. Next, we investigate the origin of the jagged likelihoods using synthetic data from the oscillator model described in eq.\ \eqref{eq:oscillator_2}.

\subsubsection{Inference for the damped, driven oscillator using an adaptive step size solver}  
We first study the effects of adaptive time step solvers on inference using the model system that was introduced at the beginning of \S{}\ref{sec:3} (eq.~\eqref{eq:oscillator_2}). Here, we set the input stimulus according to
\begin{equation} \label{eq:adaptive_ft}
F(t) = \begin{cases}1, & t<t_\text{change}, \\ f_1, & t\geq t_\text{change}. \end{cases}
\end{equation}
Thus, $f_1$ controls the strength of a pulse provided to the system at $t=t_\text{change}$.

\begin{figure}
\centering
\includegraphics[width=\textwidth]{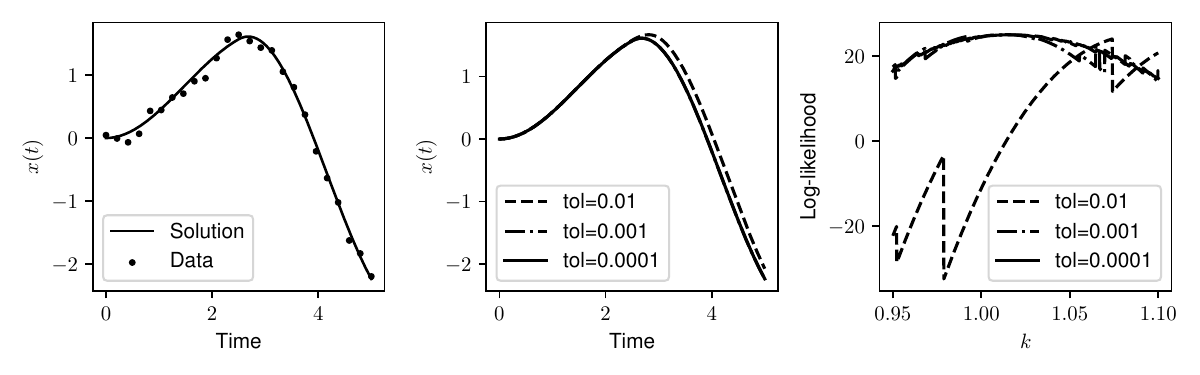}
\caption{\textbf{Damped oscillator inference using adaptive time step Runge-Kutta.} (Left) Synthetic data for the damped driven oscillator. The curved line indicates the accurate solution to the ODE with these parameters, while the points indicate the noisy data. (Center) Solution for oscillator computed using an RK5(4) solver with three different choices for the relative tolerance (indicated by $\text{tol}$ in the legend). (Right) Log-likelihood for the parameter $k$ calculated from the noisy data, with all other parameters held at their true values. The log-likelihood was calculated from eq.~\eqref{eq:iid_gaussian_ll} using an RK5(4) solver with three different choices for the tolerance.}
\label{fig:rk_oscillator_fixed_problem}
\end{figure}

First, we consider the problem where $f_1=-1$ and $t_\text{change}=2.5$ for different choices of the RK5(4) solver tolerance. $25$ evenly spaced data points were generated from and including $t=0$ to $t=5$ from the model with an accurate solver (the RK5(4) solver with relative tolerance set to $10^{-8}$), using true parameter values $k=1$, $c=0.2$, $m=1$ and an initial condition of $x(t=0)=0$, $\dot{x}(t=0)=0$. Then, IID Gaussian noise was added to the solution at each of the sampled locations with $\sigma=0.1$. Holding all other parameters fixed at their true values, the log-likelihood was calculated for a range of values of $k$, using the RK5(4) solver with various tolerances. These results are shown in Figure \ref{fig:rk_oscillator_fixed_problem}. At insufficient tolerances, the log-likelihood surface exhibits significant erroneous jaggedness. Notably, visual changes between the forward simulations are minor even at tolerances which cause drastic differences in the log-likelihood.

\begin{figure}
\centering
\includegraphics[width=\textwidth]{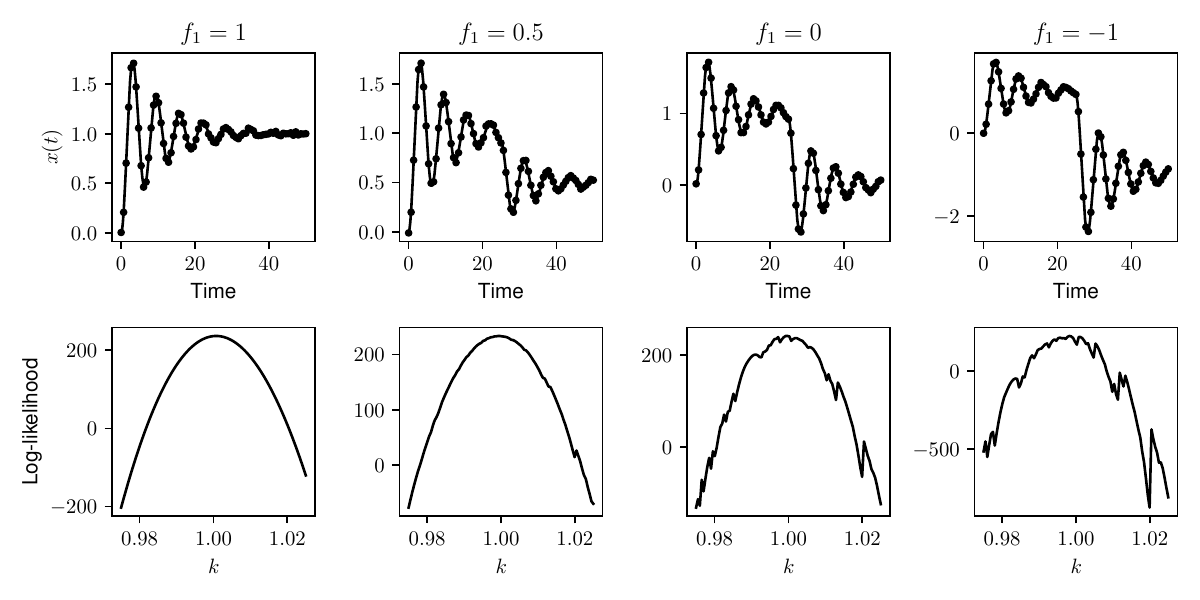}
\caption{\textbf{Damped oscillator model: forward simulations and inference using an adaptive solver.} Time series data and parameter likelihood surfaces are shown for four values of $f_1$ in the oscillator problem: eqs.\ \eqref{eq:oscillator_2} \& \eqref{eq:adaptive_ft}. For each value of $f_1$, the top plot shows the accurate ODE solution (line) and noisy synthetic data (points) generated from it. The bottom plot panels show the corresponding log-likelihood surface for $k$ over an interval centred on the true value, $k=1$, while all other parameters are held at their true values. For generating the likelihood surfaces, an RK5(4) solver was used with $\texttt{rtol}=10^{-3}$.}
\label{fig:rk_oscillator}
\end{figure}

Next, we fix the adaptive solver tolerance and study how introducing more rapid changes in the system's behavior affects the log-likelihood surface. In Figure~\ref{fig:rk_oscillator}, we fix $t_\text{change}=25$ and consider four different values of $f_1$ and plot the likelihood surface for the model parameter $m$ calculated according to an RK5(4) solver with $\texttt{rtol}=10^{-3}$. For each value of $f_1$, 75 evenly spaced data points were generated on the interval from and including $t=0$ to $t=50$, using parameter values $k=1$, $c=0.2$, and $m=1$. IID Gaussian noise was added to the solution at each of the sampled locations with $\sigma=0.01$. The likelihood was then calculated over a range of values of $k$, with all other parameters held at their correct values. For $f_1=1$, the stimulus $F(t)$ is constant over time, and the likelihood surface appears smooth. However, as $f_1$ is adjusted so the stimulus is a stronger pulse, the likelihood becomes jagged with large deviations away from the true likelihood surface. (This is an example of a challenging RHS which could be made more tractable for inference using smoothing approximations, which we analyze in \S{}S3.) Overall, these results indicate that the more rapid the changes in a system's behavior, generally the tighter solver tolerances are required to solve the inverse problem.

\begin{figure}
\centering
\includegraphics[width=.85\textwidth]{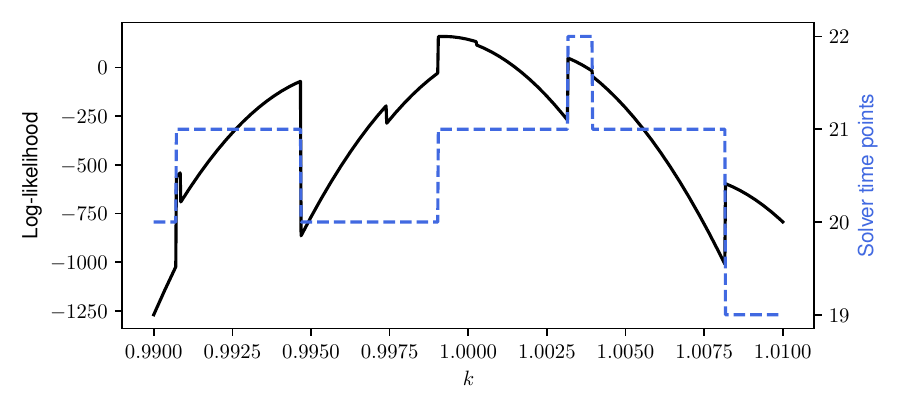}
\caption{\textbf{Damped oscillator model: likelihood discontinuities caused by variation in the number of adaptive steps.} The log-likelihood surface for the parameter $5$ in the oscillator problem (black solid line) and the number of time points used by the adaptive step size ODE solver in the calculation of each value of the likelihood (blue dashed line) are shown. An RK5(4) solver was used with $\texttt{rtol}=10^{-3}$.}
\label{fig:num_time_pts}
\end{figure}

A fundamental point to note is that these inaccuracies arise because different values of the parameters represent different forward problems, and the solver selects a different sequence of step sizes for each. When the solution contains regions of rapid change, differences in the positions of the solver time steps, and, particularly, the inevitably discontinuous jumps in the total number of time steps used by the solver, cause errors in the likelihood. This phenomenon is investigated more closely in Figure~\ref{fig:num_time_pts}. For this study, the oscillator model eq.\ \eqref{eq:oscillator_2} was again used. 
50 evenly spaced data points were generated on the time interval from and including $t=0$ to $t=10$, with $t_\text{change}=5$ and $f_1=-5$, using parameter values $m=1$, $c=0.2$, and $k=1$. IID Gaussian noise was added to the solution at each of the sampled locations with $\sigma=0.01$. The likelihood for $k$ was calculated as before and is plotted in Figure~\ref{fig:num_time_pts}. In this case, the figure is restricted to a very narrow range of $k$ values, and the total number of time points selected by the adaptive solver for the calculation of the likelihood at each value of $k$ is overlaid on the plot. Here, the large jumps in the likelihood correspond to the addition or removal of a solver time point. Smaller spikes and jaggedness where the total number of solver time points is constant correspond to shifting of the solver time points.

\subsubsection{Effect of jaggedness on inference algorithms}
The jagged spikes appearing in the likelihood surface as a result of insufficiently accurate adaptive step size solvers plague computational inference algorithms. A common approach to Bayesian inference is to use the Metropolis MCMC algorithm, or variants of it \cite{gelman2013bayesian}. This algorithm generates a sequence of parameter values via a Markov chain whose stationary distribution is the posterior distribution of the parameters. Given the most recent parameter values in the chain $\theta^\text{old}$, the Metropolis algorithm proposes new parameter values $\theta^\text{prop}$ according to a proposal distribution and then accepts $\theta^\text{prop}$ with a probability of:
$$
r = \min\left(1, \frac{p(\theta^\text{prop})}{p(\theta^\text{old})}  \frac{p(y|\theta^\text{prop})}{p(y|\theta^\text{old})} \right),
$$
where $p(\theta^\text{prop})$ is the prior and $p(y|\theta^\text{prop})$ is the likelihood. To illustrate the detrimental effects of jagged errors in the likelihood, we consider a situation where $\theta^\text{old}$ and $\theta^\text{prop}$ have identical values under the prior and the accurately computed likelihood (this is a plausible assumption when $\theta^\text{old}$ and $\theta^\text{prop}$ are nearby), but we assume that the log-likelihoods at these two parameter values computed using the numerical approximation differ by some factor $c$ driven by numerical error in the adaptive step size solver (i.e., $\log p(y|\theta^\text{prop}) = \log p(y|\theta^\text{old}) - c$, for $c>0$). This assumption of a jump in computed likelihood values at nearby parameter values is analogous to the spikes appearing in the log-likelihood in our results in Figures \ref{fig:rk_oscillator} and \ref{fig:num_time_pts}. 

Under these assumptions, $\log r = -c$ or $r = \exp(-c)$. For a value of $c=10$ (smaller than many of the magnitudes of spikes observed in our results), the probability of accepting the proposal is less than 1 in 20,000. Even a relatively small jump of magnitude $c=3$ will be traversed by the sampler with a probability of only about 5\%. Although these computations are based on simplistic assumptions, they suggest that even minor warping of the log-likelihood may severely restrict the ability of a Metropolis-Hastings sampler (or similar inference algorithm) to traverse the parameter space efficiently.

\subsection{The impact of observation error magnitude on inference and sampling performance}

\begin{figure}
\centering
\includegraphics[width=\textwidth]{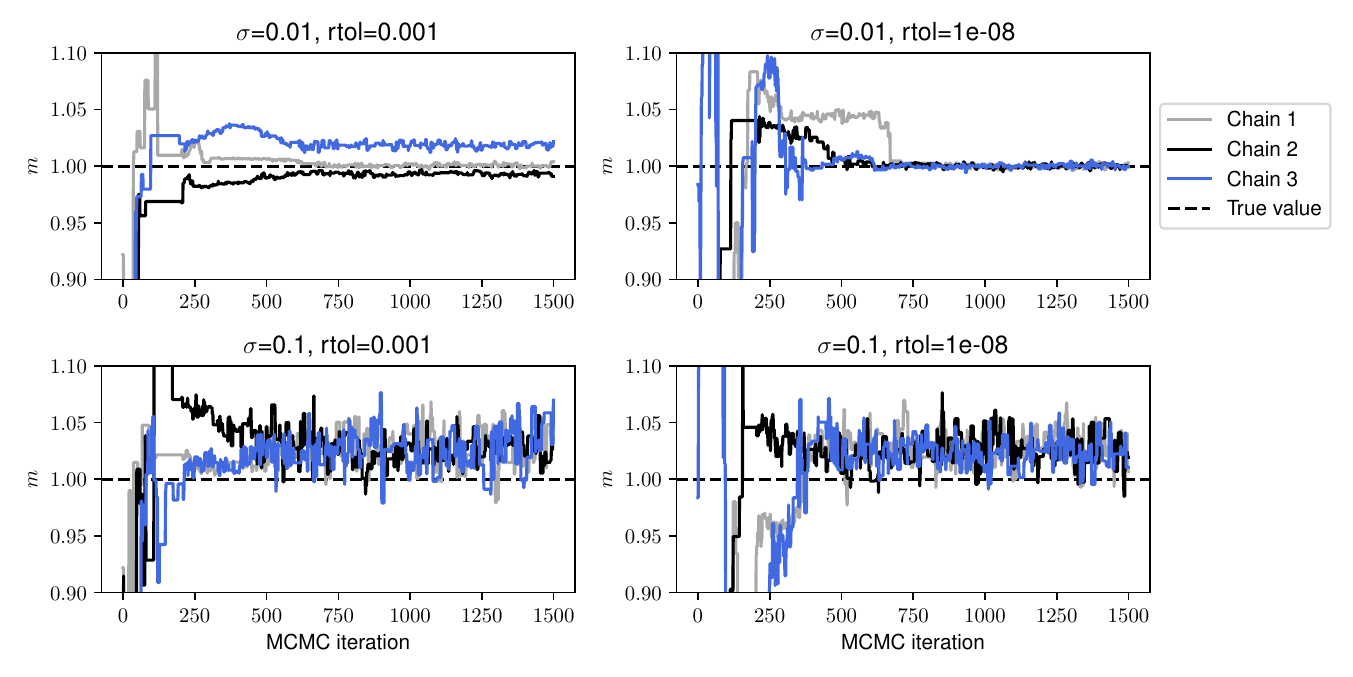}
\caption{\textbf{Effect of noise on MCMC convergence.} Data were generated according to the same specifications as for Figure~\ref{fig:rk_oscillator}, with $f_1=-1$, and the indicated values of $\sigma$. Inference was performed for the three parameters $m$, $c$, and $k$, as well as $\sigma$, via adaptive covariance MCMC~\cite{haario2001adaptive} with three independent chains initialized at random samples from the prior (uniform on $[0.1, 1.5]$ for the model parameters, and uniform on $[0, 1]$ for $\sigma$). 1500 MCMC iterations were performed. The plots show the three chains for the $m$ parameter. (Left) Forward simulation was performed using the RK5(4) solver with $\texttt{rtol}=10^{-3}$. (Right) Forward simulation was performed using the RK5(4) solver with $\texttt{rtol}=10^{-8}$.}
\label{fig:mcmc}
\end{figure}

In this section, we empirically study the effects of different levels of observation noise on inference. 
We performed Bayesian inference using MCMC for the oscillator problem with varying levels of noise in the data. We considered two values of $\sigma$ ($0.01$ and $0.1$) to generate the data, fixed $f_1=-1$, and otherwise generated data exactly as described for Figure~\ref{fig:rk_oscillator}. We set a uniform prior on $[0.1, 1.5]$ for the three model parameters $m$, $c$, and $k$, and a uniform prior on $[0, 1]$ for the $\sigma$. Three MCMC chains were run, initialized at random samples from the prior (with the same MCMC starting point being used for both choices of the true $\sigma$). 1500 iterations of MCMC were performed using the Haario-Bardenet adaptive covariance algorithm as implemented in PINTS to sample from the posterior~\cite{haario2001adaptive, clerx2018probabilistic}. The MCMC chains for the $m$ parameter are plotted in the left column of Figure~\ref{fig:mcmc} using the RK5(4) solver with $\texttt{rtol}=10^{-3}$, while the right column of Figure~\ref{fig:mcmc} shows the chains using the same solver but with more accurate tolerances of $\texttt{rtol}=10^{-8}$.

At the lowest noise level considered ($\sigma=0.01$), the three MCMC chains using the less accurate solver move towards the true value of the parameter but fail to mix with each other. Instead, each chain remains stuck in a narrow region of parameter space near the true parameter value, corresponding to the phantom local maxima in the likelihood surface observed in Figure~\ref{fig:rk_oscillator}. Reducing the solver tolerance to $10^{-8}$ was, however, sufficient to ensure chain mixing, indicating that the lack of convergence was purely an artefact of using an inaccurate solver. At the highest level of noise considered here ($\sigma=0.1$), the three MCMC chains mix well for either tolerance choice,\footnote{We note that, for this level of noise, the centers of the sampling distributions are shifted slightly away from the true parameter value because the noise limits our ability to estimate this parameter.} which can be explained by our bound given in eq.\ \eqref{eq:ll_error}: that larger $\sigma$ values lead to gentler variation in the log-likelihood surface and so easier exploration by inference algorithms.

\section{Fixed step solvers applied to an SIR change point model of the spread of COVID-19 in Germany}  \label{sir}

A widely used class of differential equation models in epidemiology are compartmental models, which divide the population into a number of compartments representing different diseased or non-diseased states and specify the rates at which individuals move from one compartment to another \cite{van2022learning}. A simple yet commonly used example is the SIR model (susceptible-infected-recovered)~\cite{weiss2013sir}. This model keeps track of the number of susceptible individuals $S(t)$ (those who can be infected with the disease), infected individuals $I(t)$ (those who are currently infectious with the disease), and recovered individuals $R(t)$ (those who have recovered from the disease and are assumed immune). Neglecting births and deaths, the model is expressed by the following system of differential equations:
\begin{equation}  \label{eq:SIR_S}
\frac{dS}{dt} = -\lambda \frac{S I}{N}
\end{equation}
\begin{equation}  \label{eq:SIR_I}
\frac{dI}{dt} = \lambda \frac{SI}{N} - \mu I
\end{equation}
\begin{equation}  \label{eq:SIR_R}
\frac{dR}{dt} = \mu I,
\end{equation}
where $\lambda>0$ is the spreading rate of the disease, $\mu>0$ is the recovery rate, and $N>0$ is the total size of the population. The system additionally requires the specification of initial conditions for each compartment ($S(0)$, $I(0)$, $R(0$)). $I(0)$ must exceed zero for an infection to spread.


The qualitative behaviour of the SIR model can be determined by the basic reproduction number, $R_0$, where
\begin{equation*}
R_0 = \frac{\lambda}{\mu}.
\end{equation*}
Assuming that $S(0) \approx N$ and $I(0) > 0$, when $R_0 > 1$, the number of infected individuals will grow, and the epidemic will eventually infect the entire population (barring a change in $\lambda$ or $\mu$); for $R_0 < 1$, the number of infected individuals will fall. Thus, fitting an SIR model to infection data, and estimating the spreading rate $\lambda$ and reproduction number $R_0$, are important steps in understanding and predicting the progression of an epidemic. 

An extension to the standard SIR model has $\lambda$ vary over time, allowing the model to capture changes in the spread of a disease caused by behavioural changes or government policy. In the aftermath of the outbreak of COVID-19 in Europe in early 2020, an SIR model allowing changes in $\lambda$ through time was used in a high profile paper which attempted to capture the impact of major public health policy interventions on COVID-19 transmission in Germany~\cite{dehning2020inferring}. The authors used the  model eqs.~\eqref{eq:SIR_S}--\eqref{eq:SIR_R}, discretised with a one day time step, equivalent to a Forward Euler solver with $\Delta t=1$:
\begin{equation}  \label{eq:S_discrete}
S_t = S_{t-1} -\lambda(t) \Delta t \frac{S_{t-1} I_{t-1}}{N}
\end{equation}
\begin{equation}
I_t = I_{t-1} +  \lambda(t) \Delta t \frac{S_{t-1}I_{t-1}}{N} - \mu \Delta t I_{t-1}
\end{equation}
\begin{equation}
R_t = R_{t-1} +  \mu \Delta t I_{t-1}.
\end{equation}
The initial condition was given by an unknown parameter $I_0 = I(0)$. The system was closed with $R(0)=0$ and $S(0) =N - I_0$. The spreading rate $\lambda$ was assumed to be a continuous function of time and was allowed to shift at three time points, whose locations were estimated from the data. Specifically, these three time points, $t_i, i \in \{1, 2, 3\}$ denoted the times at which $\lambda$ began to (linearly) change to a new, constant value, and the time taken for these shifts was dictated by durations $d_i$. The $\lambda$ profile then has the following piecewise representation:
\begin{equation*}
\lambda(t) = \begin{cases} \lambda_0, & t<t_1, \\ \lambda_0 + \frac{\lambda_1 - \lambda_0}{d_1} (t - t_1), & t_1 \leq t < t_1 + d_1, \\ \lambda_1, & t_1+d_1 \leq t < t_2, \\  \lambda_1 + \frac{\lambda_2 - \lambda_1}{d_2} (t - t_2), & t_2 \leq t < t_2 + d_2, \\ \lambda_2, & t_2+d_2 \leq t < t_3, \\ \lambda_2, + \frac{\lambda_3 - \lambda_2}{d_3} (t - t_3) & t_3 \leq t < t_3 + d_3, \\ \lambda_3, & t_3+d_3 \leq t. \end{cases}
\end{equation*}
Additional features of the model included a reporting delay and a weekly modulation. The reporting delay was characterised by a single parameter $D$ indicating the number of days between the time at which new infections occur and the time at which they are reported. The modulation accounts for the weekly periodicity evident in the data and is characterised by two parameters $f_w$ and $\Phi_w$. This significant periodicity likely arises from processes involved in the reporting of COVID-19 cases and deaths \cite{Gallagher2023.06.13.23290903}. Specifically, cases $C_t$ are modelled by:
\begin{equation}  \label{eq:delay}
C_t = (1-f(t)) I^\text{new}_{t-D},
\end{equation}
where
\begin{equation}  \label{eq:weekly}
f(t) = (1 - f_w)\left(1 - \left|\sin \left( \frac{\pi}{7} t - \frac{1}{2} \Phi_w \right) \right| \right),
\end{equation}
where $I^\text{new}_t = S_{t-1} - S_{t}$. \cite{dehning2020inferring} assumed a Student-t distribution with four degrees of freedom and multiplicative noise for the likelihood, such that the likelihood for observed cases $\hat{C_t}$ was given by:
\begin{equation*}
p(\hat{C_t} | \theta, \sigma) = \text{Student-t}_{\nu=4} (\text{mean}=C_t(\theta), \text{scale}=\sigma \sqrt{C_t(\theta)}),
\end{equation*}
where $\theta=(\lambda_0, \lambda_1, \lambda_2, \lambda_3, t_1, t_2, t_3, d_1, d_2, d_3, \mu, D, I_0, f_w, \Phi_w, \sigma)$ is the full vector of parameters for the differential equation model, and $C_t(\theta)$ is the deterministic solution which may be computed using a range of different time steps. The prior distributions for the parameters are given in Table~S1 (supplementary information).

\subsection{Effect of time step on the forward solution}

We first study the effect of assuming $\Delta t=1 \text{ day}$ on forward simulations of the model. We set up the forward simulations using the same settings that \cite{dehning2020inferring} used to generate their Figure 2. The parameters of an SIR model without change points or weekly modulation (i.e., a single value of $\lambda$, $\mu$, $D$ $I_0$, and $\sigma$) were inferred from an early period of the German daily reported COVID-19 cases, from 2 March 2020 to 15 March 2020. The posterior median values of these parameters (excepting $\lambda$) were then used to generate forward simulations according to the full model without weekly modulation (eqs.~\eqref{eq:S_discrete}--\eqref{eq:delay}), with one change point, and pre-specified values of $\lambda_0$ and $\lambda_1$.

\begin{figure}
\centering
\includegraphics[width=.95\textwidth]{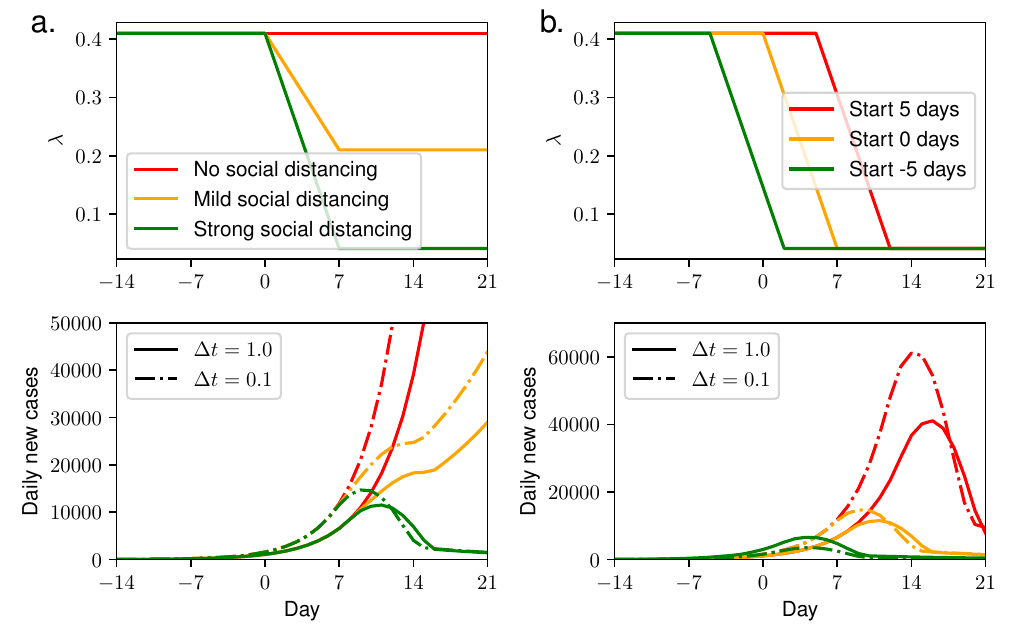}
\caption{\textbf{COVID-19 model: forward simulations using Forward Euler.} In both (a) and (b), the top panel shows three different pre-specified trajectories of $\lambda(t)$, and the bottom panel shows the number of daily cases resulting from these trajectories for each choice of the time step $\Delta t$.}
\label{fig:forecasting}
\end{figure}

As in~\cite{dehning2020inferring}, the first set of simulations considered how different levels of social restrictions could influence the course of disease transmission, as measured by cases. Three levels of social restrictions (assumed to be captured by different $\lambda$ values) are considered, which each yield two sets of simulations: one corresponding to Forward Euler with $\Delta t=1 \text{ day}$ (as in~\cite{dehning2020inferring}) and another with $\Delta t=0.1$ days. We choose $0.1$ days as the more accurate comparator method, as further refinement of the step size yields little change in forward solutions or inference results but is increasingly costly to run. The results of this are shown in Figure~\ref{fig:forecasting}A. Our second set of simulations, shown in Figure~\ref{fig:forecasting}B, considered only our ``strong'' social distancing scenario and explored three different times at which the change in $\lambda$ might occur (e.g., if a public health intervention were implemented at different times). These simulations illustrate how using a large time step generally leads to a substantial underestimation of case counts for a given choice of $\lambda(t)$, particularly during the (crucial) growth phase of the epidemic.



\begin{figure}
\centering
\includegraphics[width=\textwidth]{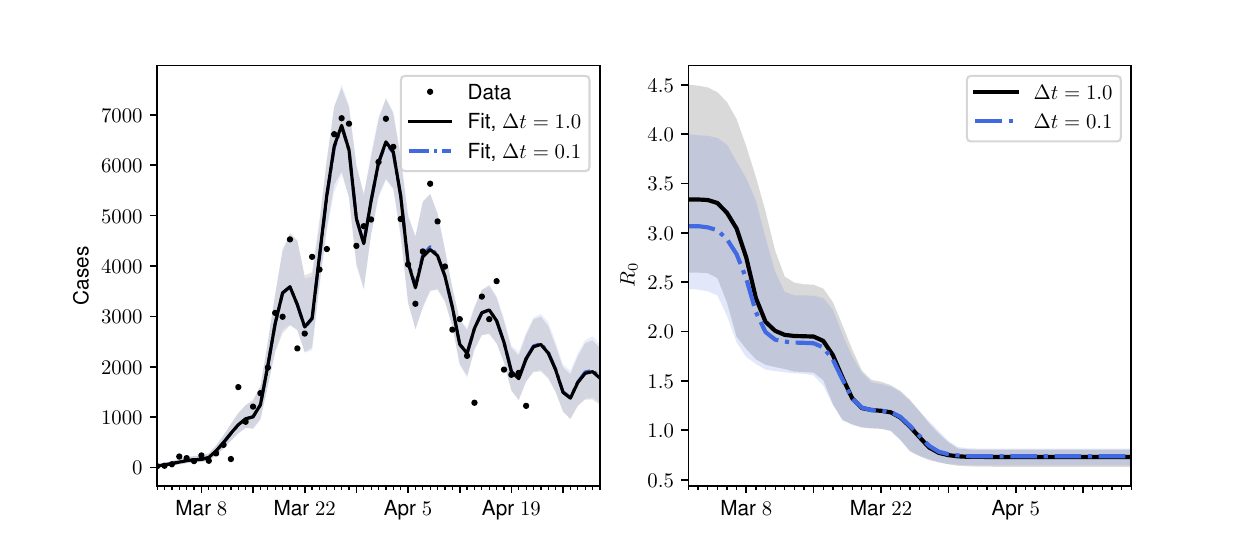}
\caption{\textbf{COVID-19 model: inference using Forward Euler} (Left) Real data and model fits for the number of daily COVID-19 cases in Germany over the period 2 March 2020 to 21 April 2020. Note that the model fits for $\Delta t=1$ and $\Delta t=0.1$ overlap almost completely. (Right) Inferred basic reproduction number over time for the Germany COVID-19 data, using the SIR model with change points in $\lambda$ (eqs.~\eqref{eq:S_discrete}--\eqref{eq:weekly}) and two different values for the ODE solver time step, $\Delta t$. In both panels, lines indicate the posterior median and shaded regions indicate the central $95\%$ of the posterior.}
\label{fig:covid_inference}
\end{figure}

\subsection{Effect of time step on the posterior distributions}
We also studied the effect of the time step on parameter inference for the full model (eqs.~\eqref{eq:S_discrete}--\eqref{eq:weekly}) using the German daily cases data from 2 March 2020 to 21 April 2020 as was done in~\cite{dehning2020inferring}. Inference was performed using the PyMC3 No-U-Turn MCMC Sampler (NUTS)~\cite{salvatier2016probabilistic,gelman2013bayesian} using the model developed by~\cite{dehning2020inferring}, modified to allow the $0.1$ day step size. To initialize the chains, automatic differentiation variational inference~\cite{kucukelbir2017automatic} as implemented in PyMC3~\cite{salvatier2016probabilistic} was performed to generate an approximate posterior (which, however, does not capture correlations between the parameters). Four MCMC chains were then initialized by sampling from this approximation of the posterior. The chains were run for $500$ iterations of NUTS, with the first $100$ discarded as burn-in, and convergence assessed by requiring that $\hat{R}<1.05$~\cite{gelman2013bayesian}. These results are shown in Figure~\ref{fig:covid_inference}.

Both models achieve a near identical visual fit to the data, using the median values of the recovered parameters. However, the parameter estimates of the two models differed. We focus on the posterior distribution for the basic reproduction number $R_0$, which is calculated using the MCMC samples of the joint posterior for $(\lambda, \mu)$. The one day time step results in overestimation of $R_0$ (by approximately 10\% relative to the $0.1$ day time step) during the early stages of the epidemic (i.e., before the first change point). This is because, during the growth phase of the epidemic, the larger time step results in slower growth for a given $\lambda$ value (cf.\ Figure~\ref{fig:forecasting}), meaning a larger $\lambda$ value is estimated to compensate. During the later stages of the epidemic, the values of $R_0$ are more similar between the two models. Additionally, the change point locations are not much affected by the choice of time step (though, this is expected as the change points have fairly informative priors).

Our results indicate that while the discrete version of the SIR change point model using $\Delta t = 1$ appears visually to obtain a good fit to German COVID-19 data, the growth parameters of the discrete model using this time step vary markedly from those recovered using $\Delta t=0.1$, and thus care should be taken in the deployment of such discrete models and the reporting of their results.

\section{Numerical errors arising in rainfall-runoff models of river streamflow data}  \label{hydrology}

In this section we use real data from the French Broad River at Asheville, North Carolina to investigate the impact of adaptive solvers in performing inference for rainfall-runoff models used in hydrology~\cite{schoups2010formal, schoups2010corruption}.

Rainfall-runoff models divide the flow of water through a river basin into several spatially grouped components representing different hydrological processes. The model we consider here is governed by a system of five ODEs:
\begin{equation}  \label{eq:si}
\frac{dS_i}{dt} = \text{Precip}(t) - \text{InterceptEvap}(t) - \text{EffectPrecip}(t)
\end{equation}
\begin{equation}
\frac{dS_u}{dt} = \text{EffectPrecip}(t) - \text{UnsatEvap}(t) - \text{Percolation}(t) - \text{Runoff}(t)
\end{equation}
\begin{equation}
\frac{dS_s}{dt} = \text{Percolation}(t) - \text{SlowStream}(t)
\end{equation}
\begin{equation}
\frac{dS_f}{dt} = \text{Runoff}(t) - \text{FastStream}(t)
\end{equation}
\begin{equation}  \label{eq:z}
\frac{dz}{dt} = \text{SlowStream}(t) + \text{FastStream}(t),
\end{equation}
Each term in this equation is defined in Table~S2 (supplementary information), and the seven unknown parameters of the model and their prior distributions are defined in Table~S3 (supplementary information). The data consist of daily streamflow measurements ($dz/dt$), and the authors assume an IID Gaussian likelihood with unknown standard deviation $\sigma$.

\begin{figure}
\centering
\includegraphics[width=\textwidth]{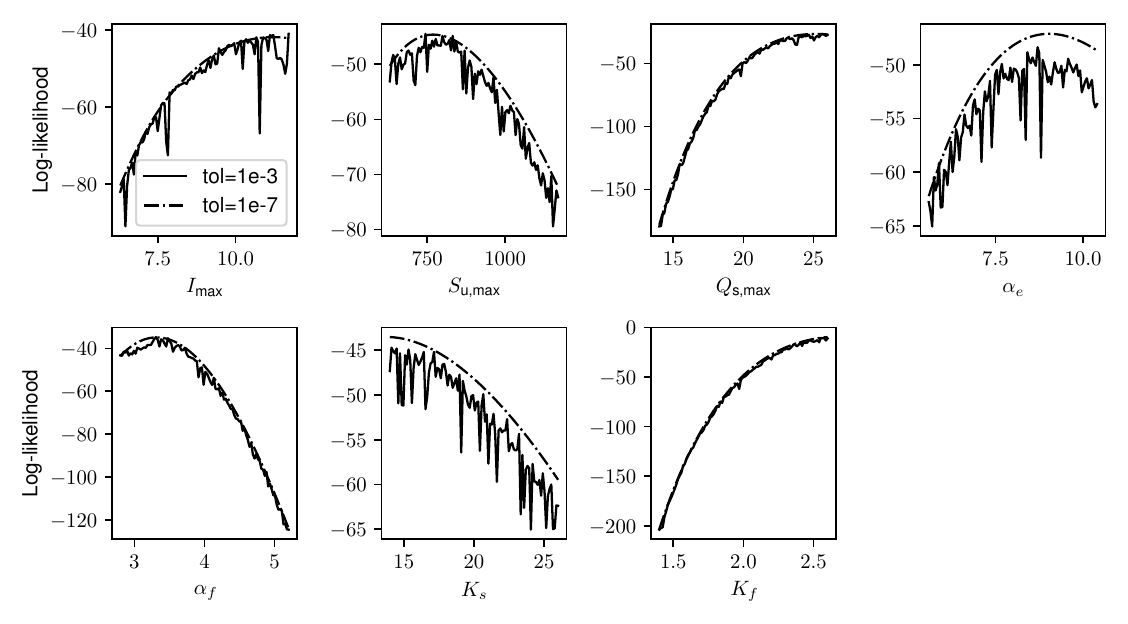}
\caption{\textbf{Rainfall run-off model: inference using accurate and inaccurate adaptive solvers.} Here, we plot the likelihood surface for each parameter for the rainfall-runoff model defined by eqs.\ \eqref{eq:si}--\eqref{eq:z}. For each parameter, the solid line indicates the likelihood calculated using an RK3(2) adaptive solver with $\texttt{rtol}=\texttt{atol}=10^{-3}$, while the dashed line indicates the likelihood calculated using the CVODE adaptive solver with $\texttt{rtol}=\texttt{atol}=10^{-7}$.}
\label{fig:rainfallrunoff}
\end{figure}

Previous work has shown that using large time steps with such hydrological models can bias inferences~\cite{kavetski2003semidistributed}. We show that using an adaptive step size method (as suggested by~\cite{schoups2010corruption}) can also cause inaccurate inference results, unless the error is tightly controlled.

Using an accurate ODE solver (the CVODE multistep solver from the SUNDIALS library~\cite{hindmarsh2005sundials} with $\texttt{rtol}=\texttt{atol}=10^{-7}$), we obtained the posterior distributions for the seven parameters of the model, using USGS data for the streamflow at Asheville, North Carolina (USGS station 03451500) over a 200 day period starting 1 January 1960. Sampling was performed using the Dream multi-chain MCMC algorithm as implemented in PINTS~\cite{vrugt2009accelerating, clerx2018probabilistic}, using 6 chains with each initialized by a sample from the prior (Table~S3, supplementary information). $25000$ MCMC iterations were performed, and convergence of the chains was assessed by requiring that $\hat{R}<1.05$~\cite{gelman2013bayesian}. In Figure~\ref{fig:rainfallrunoff}, we plot the likelihood surfaces of the parameters for slices through parameter space near the estimated posterior medians. Likelihood surfaces are plotted for two adaptive step size solvers: the RK3(2) solver from SciPy with $\texttt{rtol}=\texttt{atol}=10^{-3}$, and the CVODE solver as described above. For all parameters, the $10^{-3}$ tolerance solver causes highly jagged likelihoods, of sufficient magnitude to interfere with inference via MCMC or maximum likelihood estimation. This is in accordance with our earlier results using the oscillator model in \S{}\ref{sec:3}, as rapid changes in the RHS cause spurious jaggedness in the computed likelihood. The likelihoods calculated using the more accurate solver have similar broadscale shapes but are smooth enough for accurate inference to be performed.

\section{Discussion}



Inaccurate solution of ODEs through either fixed time step or adaptive solvers can lead to biased inferences which are generally exacerbated when there is low observation noise. For adaptive solvers, these biases may manifest through the presence of phantom jaggedness in the likelihood surface, which can wreak havoc for inference algorithms attempting to navigate the surface. They may also lead researchers to falsely conclude that a model is unidentifiable, when, for example, the chains in an MCMC run fail to mix. They may then choose to modify the model in arbitrary ways when, in fact, all that was required to render inference soluble was a reduction in solver tolerances. Tolerances which seem good enough for forward simulation are likely insufficient for solving inverse problems. For example, a default relative tolerance of $10^{-3}$ was insufficient for both the synthetic data and real data studied in \S{}\ref{adaptive} and \S{}\ref{hydrology}. When using an ODE solver library to perform inference, default settings may well not suffice and, ideally, the solver tolerance should be set by inspection of the likelihood surface.

Unless there is a bifurcation in system behavior at points in parameter space, likelihood surfaces should not have abrupt discontinuities. So, the presence of such changes may well be an artefact of using an adaptive ODE solver with insufficient tolerances. MCMC and optimisation algorithms could be augmented by monitoring for such jumps and warning the user should they occur.

ODEs involving discontinuous RHS functions are known to be particularly challenging to solve accurately. Our results indicate that RHS functions involving rapid changes over time, such as those involving discontinuities, also curse computational inference when adaptive ODE solvers are used. However, our results in \S{}S3 also indicate that errors in the likelihood arising from discontinuous RHS functions can be ameliorated through the use of simple smoothing approximations---a potentially more computationally efficient alternative to increasing tolerances. We argue that in many scientific systems such smoothing approximations are additionally more realistic descriptions of the phenomena being modelled. Although the appropriate degree of smoothing may be difficult to determine in general, for certain systems, the level of smoothing can be tuned based on knowledge of the process being modelled.

Much of the work on error control for ODE solvers has focused merely on the accuracy of the forward problem. The accuracies of widely used ODE solvers are typically tuned via their step sizes or local truncation error tolerances, but these are not the most relevant quantities for inference---instead, it is the error in the log-likelihood which must be controlled. ODE solvers which control the error on the log-likelihood directly would avoid much of thie problems highlighted in this paper, and we suggest this as a fruitful resarch direction.

\section{End section}

\subsection{Data, code and materials}

The code to perform the computer experiments presented in this paper was written in Python 3.7 and is available in an open source Python library at \url{https://github.com/rccreswell/ode_inference}. To run the COVID-19 simulations, we adapted the software library developed by \cite{dehning2020inferring}. The version of the code including our modifications is available at \url{https://github.com/rccreswell/covid19_inference_forecast}.

\subsection{Authors' contributions}

R.C.: conceptualization, formal analysis, investigation, methodology, software, validation, visualization, writing—original draft, writing—review and editing; K.M.S.: formal analysis, investigation, methodology, software, visualization, writing—original draft, writing—review and editing; B.L.: conceptualization, formal analysis, investigation, methodology, project administration, supervision, validation, writing—original draft, writing—review and editing; G.R.M: conceptualization, investigation, methodology, writing—review and editing; C.L.L.: conceptualization, investigation, methodology,  supervision, writing—original draft, writing—review and editing; S.T.: formal analysis, investigation, methodology, writing—original draft, writing—review and editing; M.R.: conceptualization, investigation, supervision, methodology, writing—original draft, writing—review and editing; D.J.G.: conceptualization, formal analysis, investigation, methodology, project administration, supervision, writing—original draft, writing—review and editing.

\subsection{Conflict of interest declaration}
We declare that we have no competing interests.

\subsection{Funding}
R.C.\ acknowledges support from the EPSRC via a doctoral training partnership studentship in the Department of Computer Science at the University of Oxford. K.M.S.\ and D.J.G.\ acknowledge funding from the EPSRC CDT in Sustainable
Approaches to Biomedical Science: Responsible and Reproducible Research - SABS:R3
(EP/S024093/1). C.L.L.\ acknowledges support from the Science and Technology Development Fund, Macao SAR (FDCT) (reference number 0048/2022/A) and support from the University of Macau via a UM Macao Fellowship.

\bibliographystyle{apalike}
\bibliography{b}

\newpage

\setcounter{section}{0}
\makeatletter 
\renewcommand{\thesection}{S\@arabic\c@section}
\makeatother

\setcounter{figure}{0}
\makeatletter 
\renewcommand{\thefigure}{S\@arabic\c@figure}
\makeatother

\setcounter{table}{0}
\makeatletter 
\renewcommand{\thetable}{S\@arabic\c@table}
\makeatother

\section{Proof of bound on the error in the log-likelihood}  \label{proof}

As stated in the main text, we denote the log-likelihood computed using the numerical approximation to the forward map by:
\begin{equation*} 
\mathcal{L}' = -\frac{N}{2} \log(2\pi) - \frac{N}{2} \log(\sigma^2) - \frac{1}{2\sigma^2} \sum_{i=1}^N (y_i - g(\hat{x}_i))^2,
\end{equation*}
while the log-likelihood computed using the hypothetical, true solution to the ODE is:
\begin{equation*} 
\mathcal{L} = -\frac{N}{2} \log(2\pi) - \frac{N}{2} \log(\sigma^2) - \frac{1}{2\sigma^2} \sum_{i=1}^N (y_i - g(x(t_i;\theta)))^2,
\end{equation*}
where $x(t; \theta)$ is the true solution to the ODE, $\{\hat{x}_i\}_{i=1}^N$ is the approximate solution to the ODE at the data time points $\{t_i\}_{i=1}^N$, and $g$ is the observation operator.

For brevity, let $g_i = g(x(t_i; \theta))$ and $\hat{g}_i = g(\hat{x}_i)$. We have:
\begin{align*}  
    \mathcal{L}-\mathcal{L}' &=  - \frac{1}{2\sigma^2} \sum_{i=1}^N (y_i - g_i)^2 + \frac{1}{2\sigma^2} \sum_{i=1}^N (y_i - \hat{g}_i)^2 
    \\
    &=  \frac{1}{2\sigma^2} \sum_{i=1}^N \left[  -g_i^2 + \hat{g}_i^2 + 2y_i (  g_i - \hat{g}_i )\right] \\
    &= \frac{1}{2\sigma^2} \sum_{i=1}^N  (  g_i - \hat{g}_i ) \left(2 (y_i - g_i) + (g_i - \hat{g}_i) \right) \\
    &= \frac{1}{2\sigma^2} \sum_{i=1}^N  \left(    (  g_i - \hat{g}_i )^2 + 2 (  g_i - \hat{g}_i ) (y_i - g_i) \right ). \numberthis \label{eq:diff1} 
\end{align*}
Taking the absolute value,
\begin{align*}
    |\mathcal{L}-\mathcal{L}'| &= \left|  \frac{1}{2\sigma^2} \sum_{i=1}^N  \left(    (  g_i - \hat{g}_i )^2 + 2 (  g_i - \hat{g}_i ) (y_i - g_i) \right )  \right|  \\
    &\leq \frac{1}{2\sigma^2}  \sum_{i=1}^N \left|  (  g_i - \hat{g}_i )^2 + 2 (  g_i - \hat{g}_i ) (y_i - g_i) \right| \\
    &\leq \frac{1}{2\sigma^2}  \sum_{i=1}^N   \left( (  g_i - \hat{g}_i )^2 + 2 \left| (  g_i - \hat{g}_i ) (y_i - g_i) \right| \right)\numberthis \label{eq:ksproof3} \\
\end{align*}
To proceed further, we impose the assumption of Lipschitz continuity of the observation function $g$ with Lipschitz constant $K$, i.e., $|g(x_1) - g(x_2)| \leq K |x_1 - x_2|$ for all $x_1, x_2 \in \mathbb{R}^l$. We thus bound:
\begin{align*}
    |g_i - \hat{g}_i| = |g(x(t_i; \theta)) - g(\hat{x_i})| &\leq K |x(t_i; \theta) - \hat{x_i} | \\
    &= K|e(t_i)|. \numberthis \label{eqn}
\end{align*}
Using this in eq.\ \eqref{eq:ksproof3},
\begin{align*}
    |\mathcal{L}-\mathcal{L}'| &\leq \frac{1}{2\sigma^2}   \sum_{i=1}^N \left( K^2 |e(t_i)|^2 + 2 K |e(t_i)| |y_i - g_i|  \right).
\end{align*}
As discussed further in the main text, this bound indicates an inverse relationship between $\sigma$ and $|\mathcal{L}-\mathcal{L}'|$ when $e(t_i)$ is held constant.

\section{Distribution of the error in the log-likelihood}  \label{proofKSattempt}
In this section, rather than deriving a bound on the \emph{absolute value} of the error in the log-likelihood arising from numerical errors in the forward solution, we study the \emph{distribution} of the difference between the true and numerically approximated log-likelihoods.

We assume that the $y_i$ are distributed according to the specified IID Gaussian likelihood at the same parameter values at which the likelihood is being evaluated, so we have
$y_i \sim N(g(x(t_i; \theta)), \sigma)$. 
For brevity, as before, let $g_i = g(x(t_i; \theta))$ and $\hat{g}_i = g(\hat{x}_i)$, and define $a_i = g_i - \hat{g}_i$. Then, we can write $y_i = g_i + \epsilon_i$, where $\epsilon_i \sim N(0, \sigma)$. Starting from \S{}S1, eq.\ \eqref{eq:diff1}, we have:
\begin{align*}  
    \mathcal{L}-\mathcal{L}'&= \frac{1}{2\sigma^2} \sum_{i=1}^N  \left(    a_i^2 + 2 a_i (y_i - g_i) \right )  \\
    &= \frac{1}{2\sigma^2} \sum_{i=1}^N  \left(    a_i^2 + 2 a_i \epsilon_i  \right )  \text{ using }  y_i = g_i + \epsilon_i \\
    &=  \frac{1}{\sigma^2} \sum_{i=1}^N a_i \epsilon_i + \frac{1}{2\sigma^2}\sum_{i=1}^N a_i^2. \numberthis \label{eq:dists1} 
\end{align*}
Noting that the first term in eq.\ \eqref{eq:dists1} is the sum of independent random variables each with mean $0$ and variance $(a_i^2 / \sigma^2)$, we obtain the result:
\begin{align*}
     \mathcal{L}-\mathcal{L}' \sim  N \left(\frac{\sum_{i=1}^N  a_i^2  }{2\sigma^2},  \frac{\sqrt{ \sum_{i=1}^N a_i^2} }{\sigma} \right) .  
\end{align*}

We have $\mathbb{E}[\mathcal{L}-\mathcal{L}'] = \frac{\sum_{i=1}^N  a_i^2}{2\sigma^2}$. Therefore, on average the numerical approximation of the log-likelihood underestimates the true log-likelihood (when both are computed at the same parameter values which generated the data). The average size of this underestimation is greater when $\sigma$ is smaller. Also, in the case that $g$ is the identity function, $\sum_{i=1}^N a_i^2 = \sum_{i=1}^N |e(t_i)|^2$, and so the average size of the underestimation decreases as the global truncation error decreases.

\section{Smoothing forcing terms to reduce numerical errors in the likelihood}  \label{sec:smoothing}
As indicated by our results in Figures 4 and 5, discontinuities in the right-hand side (RHS) of an ODE can lead to substantial errors in the likelihood when adaptive step size solvers are used. In general, errors in the likelihood  arising from numerical errors in the solution can be reduced by refining the tolerance of the adaptive solver. However, when the RHS suffers from a discontinuity, the required solver tolerance to obtain an acceptable likelihood surface may employ a prohibitively large grid of solver points. Several approaches to remove discontinuities from the RHS have been developed to enable more accurate forward simulations, including smoothing approximations and solving the ODE separately within regions where the RHS is continuous~\cite{stewart2011dynamics}. These techniques may be particularly advantageous when performing inference. In this section, we study the effects on the computation of the log-likelihood of one of these approaches, which is to smooth discontinuities in the RHS of the ODE. Smoothing is often a particularly appropriate assumption for biological models, where a continuous rather than an instant change may in fact more realistically represent the true behaviour of the system. For example, in epidemiology, interventions (such as the introduction of a vaccination campaign) may be naively represented by discontinuous step functions in the RHS of a compartmental epidemiological model; however, a function smoothly moving between two values (corresponding to the intervention reaching its full effect gradually over an appropriate period of time) is both more realistic and more tractable for numerical solvers for the forward problem~\cite{van2022learning}.

The hyberbolic tangent function ($\tanh$) is a useful smooth approximation to a step function. In the forced oscillator problem, we can use $\tanh$ to approximate the step function stimulus, eq.\ (9) (main text), with $f_1=-1$ according to:
\begin{equation} \label{eq:tanh}
    F^\text{smooth}(t) = -\tanh\left(\frac{ t- t_\text{change}}{a} \right)
\end{equation}
where $a$ is a tuning parameter controlling the level of smoothing, with larger values of $a$ leading to a more gradual change in the stimulus, and $t_\text{change}$ is the time when the stimulus changes in value.

\begin{figure}[h!]
\centering
\includegraphics[width=\textwidth]{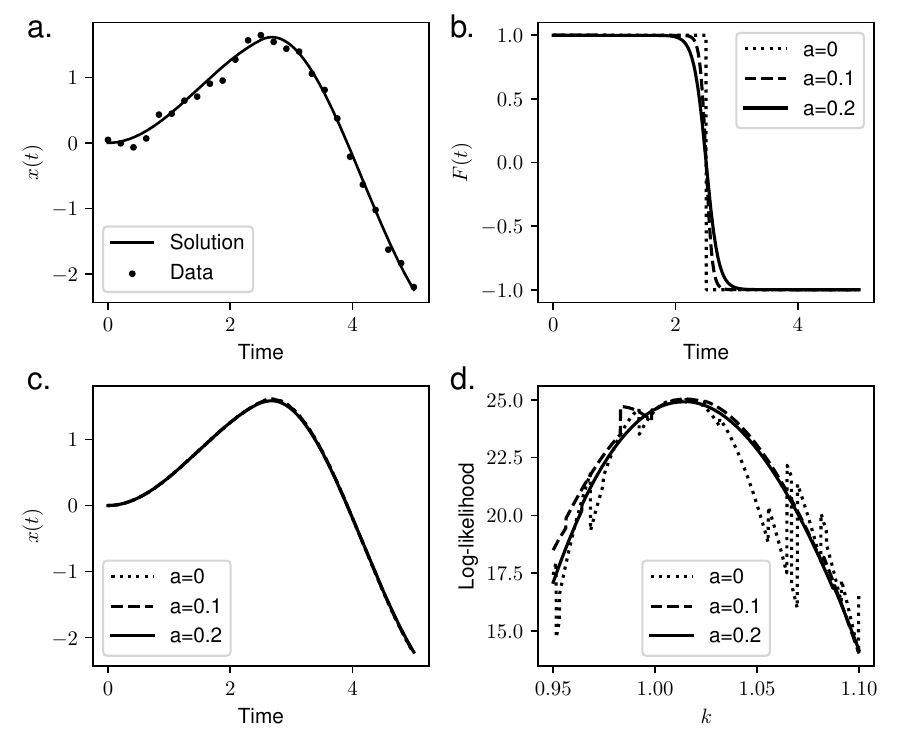}
\caption{\textbf{Effect of tanh-smoothing on likelihood surface.} (a.) Synthetic data for the damped driven oscillator. The curved line indicates the accurate solution to the ODE with these parameters, while the points indicate the noisy data. (b.) The three considered forms of the stimulus. $a=0$ indicates the unsmoothed stimulus (eq.\ (9), main text), while the positive values of $a$ indicate the tanh-smoothed stimulus according to eq.\ \eqref{eq:tanh}. (c.) Solution for oscillator computed using an RK5(4) solver with relative tolerance $10^{-3}$, with three different forms of the stimulus, at the true parameter values. (d.) Log-likelihood for the parameter $k$ calculated from the noisy data, with all other parameters held at their true values. The log-likelihood was calculated from eq.~(5) (main text) using an RK5(4) solver with relative tolerance $10^{-3}$.} 

\label{fig:smoothing}
\end{figure}

To examine the effect of the smoothing approximation on inference, we computed the likelihood surface for the $k$ parameter in the forced oscillator model using a variety of choices for the smoothing parameter, with results shown in Figure \ref{fig:smoothing}. Using $f_1=-1$ and $t_\text{change}=2.5$, $25$ evenly spaced data points were generated from and including $t=0$ to $t=5$ from the model with an accurate solver (the RK5(4) solver with relative tolerance set to $10^{-8}$), using true parameter values $k=1$, $c=0.2$, $m=1$ and an initial condition of $x(t=0)=0$, $\dot{x}(t=0)=0$. Then, IID Gaussian noise was added to the solution at each of the sampled locations with $\sigma=0.1$. Holding all other parameters fixed at their true values, the log-likelihood was calculated for a range of values of $k$, using the RK5(4) solver with relative tolerance tuned to $10^{-3}$. The likelihood was computed using both the original step function stimulus eq.\ (9) (main text) (indicated in Figure \ref{fig:smoothing} by $a=0$), as well as the smooth approximation eq. \eqref{eq:tanh} with two different choices of $a>0$. Without smoothing, we observe significant jagged biases in the likelihood, as expected due to the insufficient solver tolerance. However, with smoothing, a smooth, tractable likelihood surface is obtained despite the mediocre solver tolerance. This is despite the fact that all forward simulations are visually very similar. This is in accordance with our results in Figure~3 (main text), where even visibly small changes in the forward solution may hide the fact that there lurks substantial distortions of the likelihood surface.

\clearpage
\section{Supplementary information}





\begin{table}[h!]
\centering
\begin{tabular}{|l|l|} 
 \hline
 Parameter &  Prior \\ 
 \hline
    $\lambda_0$ & $\text{log normal}(\log(0.4), 0.5)$  \\
    $ \lambda_1$ & $\text{log normal}(\log(0.2), 0.5)$ \\
    $\lambda_2 $ & $\text{log normal}(\log(0.125), 0.5)$ \\
    $\lambda_3$ & $\text{log normal}(\log(0.0625), 0.5)$ \\
    $t_1$ & $N(2020 \text{ March } 9, 3 \text{ days})$ \\
    $t_2$ & $N(2020 \text{ March } 16, 1 \text{ day})$ \\
    $t_3$ & $N(2020 \text{ March } 23, 1 \text{ day})$ \\    
    $d_i$ & $\text{log normal}(\log(3), 0.3)$ \\
    $\mu$ & $\text{log normal}(\log(0.0625), 0.2)$ \\
    $D$ & $\text{log normal}(\log(8), 0.2)$ \\
    $I_0$ & $\text{half Cauchy}(100)$ \\
    $f_w$ & $\text{beta}(0.7, 0.17)$ \\
    $\Phi_w$ & $\text{Von-Mises}(0, 0.01)$ \\
    $\sigma$ & $\text{half Cauchy}(10)$ \\ \hline
\end{tabular}
\caption{Prior distributions for parameters in the SIR changepoint model.}
\label{table:sir_priors}
\end{table}

\clearpage
\begin{table}
\centering
\begin{tabularx}{1.1\textwidth}{|l|X|X|} 
 \hline
 Term & Definition & Description \\ 
 \hline
  $S_i$ & Interception storage & Water which strikes vegetative surfaces. \\
 \hline
  $S_u$ & Unsaturated storage & Storage of water in the soil above the water table. \\
 \hline
  $S_s$ & Slow reservoir & Water moving to the river via percolation. \\
 \hline
  $S_f$ & Fast reservoir & Water moving to the river via surface runoff. \\
 \hline
  $z$ & River discharge & Water flowed out of the river at the measuring location. \\
 \hline
 $f(S,a)$ & $\frac{1 - e^{-a  S}}{1 - e^{-a}}$ & Nonlinear flux function. \\
 \hline
 $\text{Precip}(t)$ & Precipitation & Areal precipitation in the river basin, provided as input to the model. \\
 \hline
  $\text{Evap}(t)$ & Evaporation & Evaporation from the river basin, provided as input to the model. \\
 \hline
  $\text{InterceptEvap}(t)$ & $\text{Evap}(t) f(S_i / I_\text{max}, \alpha_i)$ & Evaporation from interception. \\
  \hline
   $\text{EffectPrecip}(t)$ & $\text{Precip}(t)  f(S_i / I_\text{max}, -\alpha_i)$ & Effective precipitation which reaches unsaturated storage. \\
 \hline 
   $\text{UnsatEvap}(t)$ & $\max(0, \text{Evap}(t) - \text{InterceptEvap}(t))
                       f(S_u / S_\text{u,max}, \alpha_e)$ & Evaporation from unsaturated storage. \\
 \hline  
   $\text{Percolation}(t)$ & $Q_\text{s,max} f(S_u / S_\text{u,max}, \alpha_s)$ & Trickling of water through the ground. \\
 \hline 
   $\text{Runoff}(t)$ & $\text{EffectPrecip}(t) f(S_u / S_\text{u,max}, \alpha_f)$ & Flow of water on the surface. \\
 \hline 
   $\text{SlowStream}(t)$ & $S_s / K_s$ & Slow component of the river flow. \\
 \hline 
   $\text{FastStream}(t)$ & $S_f / K_f$ & Fast component of the river flow. \\
 \hline 
\end{tabularx}
\caption{Description of the terms which appear in the rainfall-runoff model.}
\label{table:1}
\end{table}

\clearpage
\begin{table}
\centering
\begin{tabularx}{\textwidth}{|X|X|X|} 
 \hline
 Parameter & Definition & Prior \\ 
 \hline
  $I_\text{max}$ & Maximum interception & $\text{Uniform}(0, 10)$\\
 \hline
  $S_\text{u,max}$ & Unsaturated storage capacity & $\text{Uniform}(10, 1000)$\\
 \hline
  $Q_\text{s,max}$ & Maximum percolation & $\text{Uniform}(0, 100)$\\
 \hline
  $\alpha_e$ & Evaporation flux shape & $\text{Uniform}(0, 100)$\\
  \hline
  $\alpha_f$ & Runoff flux shape  & $\text{Uniform}(-10, 10)$\\
  \hline
   $K_s$ & Slow reservoir time constant & $\text{Uniform}(0, 150)$\\
   \hline
   $K_F$ & Fast reservoir time constant & $\text{Uniform}(0, 10)$\\
 \hline
   $\alpha_s=0$ & Percolation flux shape & - \\
 \hline
   $\alpha_i=50$ & Interception flux shape & -\\
 \hline
 $\sigma$ & Noise standard deviation & \text{Uniform}(0, 10) \\ \hline
\end{tabularx}
\caption{Description of the seven unknown parameters of the model, and the two parameters with fixed values.}
\label{table:2}
\end{table}



\end{document}